\documentclass[11pt]{article}

\begin{document}

\newcommand{\ot}{\otimes}
\newcommand{\oz}{\omega}
\newcommand{\ra}{\rightarrow}
\newcommand{\op}{\oplus}
\newcommand{\lz}{\lambda}
\newcommand{\fz}{\varphi}
\newcommand{\ez}{\varepsilon}
\newcommand{\Oz}{\Omega}
\newcommand{\om}{\omega_{\mu}}
\newcommand{\on}{\omega_{\nu}}
\newcommand{\Gz}{\Gamma}
\newcommand{\gz}{\gamma}
\newcommand{\az}{\alpha}
\newcommand{\sm}{s_{\mu}}
\newcommand{\sn}{s_{\nu}}
\newcommand{\Lz}{\Lambda}
\newcommand{\sz}{\theta}
\newcommand{\pz}{\psi}
\newcommand{\Fz}{\Phi}
\newcommand{\dz}{\delta}
\newcommand{\ca}{{\cal K}}
\newcommand{\cb}{{\cal B}}
\title{\bf $AF$ Embedding of Crossed Products of Certain Graph $C^*$-Algebras by Quasi-free Actions}
  % Enter your title between curly braces
\author{ XIAOCHUN FANG }% Enter your name between curly braces

\date{}

\maketitle
\begin{abstract}
We introduce the labelling map and the quasi-free action of a
locally compact abelian group on a graph $C^*$-algebra of a
row-finite directed graph. Some necessary conditions for embedding
the crossed product to an $AF$ algebra are discussed, and one
sufficient condition is proved that if the row-finite directed
graph is constructed by possibly attaching some 1-loops to a
row-finite directed graph whose each weak connected component
is a rooted (possibly infinite) directed tree, and the labelling map is almost proper, which is proved to be a reasonable generalization of the earlier case, then the crossed product can be embedded to an $AF$ algebra.
\end{abstract}

\section*{1. Introduction}

\hspace{0.6cm} For more than ten years many important progresses
have been made in the classification of amenable $C^*$-algebras,
relative to which much attention are paid to some special
$C^*$-algebra classes, for example $AH$ algebras, purely infinite
$C^*$-algebras, crossed products of some $C^*$-dynamical systems,
quasi-diagonal $C^*$-algebras and so on. It is well known that
$AF$ embedding implies quasi-diagonality. Since M. Pimsner and D.
Voiculescu's $AF$ embedding result of irrational rotation
$C^*$-algebras (see \cite{[PV]}), there are much effort to be made
to embed the crossed products of some special $C^*$-algebras to
$AF$ algebras (for example see \cite{[Br]},
\cite{[Pi1]}--\cite{[Pu]}), which sometimes also deduce some
$K$-theory information. These special $C^*$-algebras are generally
the finite $C^*$-algebras, which is clearly a necessary condition
if the action group is discrete. Recently T. Katsura in
\cite{[Ka]} embeds successfully certain crossed products of the
Cuntz algebras ${\cal O}_n$, which are purely infinite simple
$C^*$-algebras, by quasi-free actions of a locally compact
continuous abelian group, where the quasi-free actions on ${\cal
O}_n$ have been studied for many years specially by A. Kishimoto.
Graph $C^*$-algebras are the generalization of the
 Cuntz algebras and the Cuntz-Krieger algebras. Many results about
their properties, ideal structures and $K$-theories have been
gotten for several years recently (for example see \cite{[BPRS]},
\cite{[EL]}, \cite{[F3]}, \cite{[KPR]}, \cite{[KPRR]}).

In this paper, we introduce the quasi-free action of a locally
compact abelian group $G$ with the dual $\Gamma$ on a graph
$C^*$-algebra $C^*(E)$ of a row-finite directed graph $E$, which
is defined by a labelling map $\oz$ from $E^*$ to $\Gamma$.
Keeping the embedding frame in \cite{[Ka]}, we will prove that the
crossed product of $C^*(E)$ by $G$ can be embedded to an $AF$
algebra for some special $E$ and $\oz$, which generalizes the main
result Theorem 3.8 in \cite{[Ka]} from ${\cal O}_n$ to a much
bigger $C^*$-algebra class $C^*(E)$, which contain both some
simple $C^*$-algebras and some non-simple $C^*$-algebras, and also
contain both some purely infinite $C^*$-algebras and some finite
$C^*$-algebras.

This paper is organized as follows. In section 2, we introduce the
labelling map $\oz: E^*\to \Gamma$ and the quasi-free action of
$G$ on $C^*(E)$ defined by $\oz$. From the discussion of the
necessary conditions for embedding the crossed products of
$C^*(E)$ by $G$ to $AF$ algebras,  we introduce the concept that
$\oz$ is almost proper, which is proved to be a reasonable
generalization of the corresponding concept in \cite{[Ka]}. In
section 3, we introduce the row-finite directed graphs $E$ which
are constructed by possibly attaching some 1-loops to a row-finite
directed graph whose each weak connected component is a rooted
(possibly infinite) directed tree. Then it is proved that for
these $E$, the crossed products can be embedded to $AF$ algebras,
if the corresponding labelling maps are almost proper.

\section*{2. Quasi-free actions and almost proper maps}

\hspace{0.6cm} A directed graph $E=(E^{0},E^{1},r,s)$ consists of
countable (possibly infinite) sets $E^{0}$ of vertices, $E^{1}$ of
edges, and maps $r,s: E^{1} \longrightarrow E^{0}$ identifying the
range and source of each edge. The graph is  called row-finite if
each vertex emits at most finitely many edges. We write $E^{n}$
for
 the set
of paths $\mu=e_{1}e_{2}\cdots e_{n}$ with length $|\mu|=n$, which
are sequences of edges $e_{i}$ such that $r(e_{i})=s(e_{i+1})$ for
$1\le i<n$. Then the maps $r,s$ extend naturally to
$E^{*}=\cup_{n\ge 0}E^{n}$ and $s$ extends naturally to the set of
infinite paths $\mu=e_{1}e_{2}\cdots$. In particular, we have
$r(v)=s(v)=v$ for $v\in E^{0}$. A path $\mu$ with $|\mu|\ge 1$ is
called a loop if $s(\mu)=r(\mu)$, and a loop $\mu$ is called a
1-loop if $|\mu|=1$. A vertex $v\in E^{0}$ which emits no edges is
called a sink. the relation $ \leq_E $ on $E^{0}$ is defined by $v
\leq_E w$ if and only if there is a path $\mu \in  E^{*}$ with
$s(\mu)=w$ and $r(\mu) =v$. A directed graph $E=(E^{0},E^{1},r,s)$
is called finite if both $E^0$ and $E^1$ are finite sets.

For a directed graph $E=(E^{0},E^{1},r,s)$, the weak connected
relation $\sim$ in $E^0$ is defined as following: for $v,w\in
E^0$, $v\sim w$ if and only if $v=w$ or there are $e_1, e_2,
\cdots, e_n$ in $E^1$, and $v_1, v_2, \cdots, v_n$ in $E^0$ such
that $v_0=v$, $v_n=w$, and
$\{v_{i-1},v_i\}=\{r(e_i),s(e_i)\}$($i=1,2,\cdots,n$). Clearly
$\sim$ is an equivalence relation in $E^0$. A directed graph
$F=(F^{0},F^{1},r_F,s_F)$ is called a weak connected component of
$E$ if $F^0$ is an equivalence class of $\sim$ in $E^0$,
$F^1=\{e\in E^1 |r(e)\in F^0\}=\{e\in E^1 |s(e)\in F^0\}$, and
$r_F=r|_{F^1}$, $s_F=s|_{F^1}$. A directed graph
$E=(E^{0},E^{1},r,s)$ is called a rooted directed tree if there is
a $v_0\in E^0$ with the property that there exists a unique path
in $E^*$ from $v_0$ to every other vertex in $E^0$, but no path
from $v_0$ to $v_0$.

Let $E$ be a row-finite directed graph, and let $A$ be a $C^{*}$
-algebra. A Cuntz-Krieger $E$-family in $A$ consists of a set
$\{p_{v}: v \in E^{0}\}$ of mutually orthogonal projections in $A$
and a set $\{s_{e}: e \in E^{1}\}$ of partial isometries in $A$
satisfying that $s_{e}^{*}s_{e}=p_{r(e)}$ for $e \in E^{1}$ and $
p_{v}=\sum_ {\{e:s(e)=v\}} s_{e}s_{e}^{*}$, whenever $ v$ is not a
sink.  Clearly if each $s_e$ ($e\in E^1$) is not zero, then the
product $s_{\mu}=s_{\mu_{1}}s_{\mu_{2}} \cdots s_{\mu_{n}}$, where
$\mu_i\in E^* (1\le i\le n)$ and $\mu=\mu_{1}\mu_{2}\cdots
\mu_{n}$, is non-zero
 precisely when
$\mu$ is a path in $E^{*}$. Since the range projections
${s}_{e}{s}_{e}^{*}(e\in E^1)$ are mutually orthogonal, we have
${s}_{e}^{*}{s}_{f}=0$ unless $e=f$ for any $e,f\in E^1$. For
convenience, since vertices are paths of length 0, we let
$s_{v}=p_{v}$ for $v \in E^{0}$.

 Let $E$ be a row-finite directed graph, as it is showed
  in \cite{[BPRS]},
there is a $C^{*}$-algebra $C^{*}(E)$ (called the graph
$C^*$-algebra of $E$) which is generated by a Cuntz-Krieger
$E$-family $\{s_{e},p_{v}\}$ in $C^{*}(E)$ of non-zero elements
such that, for any Cuntz-Krieger $E$-family $\{S_{e},P_{v}\}$ in
$B({\cal K})$ for some Hilbert space ${\cal K}$, there is a
representation $\pi=\pi_{S,P}$ of $A$ on ${\cal K}$ satisfying
that $\pi(s_{e})=S_{e}$, $\pi(p_{v})=P_{v}$, for all $e \in
E^{1}$, $v \in E^{0}$.
  With the convenience that
$p_{v}=s_v = s_{v}s_{v}^{*}$ for $v \in E^{0}$, we have
$C^{*}(E)=\overline{\rm span}\{s_{\mu}s_{\nu}^{*}: \ \mu,\nu \in
E^{*} \;{\rm and}\; r(\mu)=r(\nu)\in E^0\}$.\\

 {\sl{\bf Definition 2.1.} Let $G$ be a (we always assume second countable)
 locally compact abelian group with the dual $\Gz$, which is also
  a locally compact abelian group, a map $\oz: E^*=\cup_{n\ge 0}E^n\to
  \Gamma$ is called a labelling map, if
$\oz(\mu)=\oz(e_1)+\oz(e_2)+\cdots +\oz(e_n)$ for $\mu=e_1 e_2
\cdots e_n\in E^n$ and
  $\oz(\mu)=1_{\Gamma}$ for $\mu\in E^0$, where $1_{\Gamma}$ is
  the unit of $\Gamma$. }

It is clear that the labelling map $\oz$ is determined by
  $\oz|_{E^1}$, which is really just labelling of the edges of the directed graph $E$ by
  elements of the group $\Gamma$. For convenience, we denote $\oz(\mu)$ by $\oz_{\mu}$,
   and clearly the image $\oz(E^*)$
  of $\oz$ is a countable set.\\

  For any $t\in G$, let $\widetilde{s_e}=(t,\oz_e)s_e$,
  $\widetilde{p_v}=p_v$, it is easy to see that
  $(\widetilde{s_e},\widetilde{p_v})$ is a Cuntz-Krieger E-family in
  $C^*(E)$. We have then an endmorphism
  $\az^{\oz}_t : C^*(E)\to C^*(E)$ with $\az^{\oz}_t(s_e)=\widetilde{s_e}$ and
  $\az^{\oz}_t(p_v)=\widetilde{p_v}$. Since
  $\az^{\oz}_{-t}=(\az^{\oz}_t)^{-1}$, $\az^{\oz}_t$ is an
  automorphism of $C^*(E)$, and moreover $(C^*(E), G, \az^{\oz})$
  is a $C^*$-dynamical system. It is easy to see that
$\az_t^{\oz}(s_{\mu}s_{\nu}^*)=(t,\oz_{\mu}-\oz_{\nu})s_{\mu}s_{\nu}^*$
for any $\mu,\nu\in E^*$. We call the action $\az^{\oz}$ of $G$ on
$C^*(E)$ the quasi-free action. Clearly the gauge action is a
special quasi-free action.

  Viewing $C^*(E)$ as a $C^*$-algebra on a Hilbert space ${\cal H}$,
since $G$ is abelian, the regular representation
  $\tilde{\pi}\times\lz$ of $C^*(E)\times_{\az^{\oz}}G$
 on $L^2(G,{\cal H})$ is faithful. Let
   $u: L^2(G,{\cal H})\to L^2(\Gamma,{\cal H})$ be the Fourier
   transformation, $(u\xi)(\sigma)=\int_G(t,\sigma)\xi(t)dt$,
   $\forall \sigma\in \Gamma$, $\xi\in L^2(G,{\cal H})$, then for
   any $\mu\in E^*$, $g\in C_c(G)\subseteq C^*(G)$,
   $\eta\in L^2(\Gamma,{\cal H})$, we have
  $(u(\tilde{\pi}\times\lz)(g)u^*\eta)(\sigma)={\hat
   g}(\sigma)\eta(\sigma)$, and
   $(u(\tilde{\pi}\times\lz)(s_{\mu})u^*\eta)(\sigma)
   =s_{\mu}(\eta(\sigma-\oz_{\mu}))$, where ${\hat
   g}$ is the Fourier
   transformation of $g$, ${\hat
   g}(\sigma)=\int_G(t,\sigma)g(t)dt$. Under the spatial
   isomorphism defined by $u$ above, $C^*(E)\times_{\az^{\oz}}G$
   then can be considered as a
   $C^*$-algebra on $L^2(\Gamma, {\cal H})$. For the convenience,
   through the rest of the paper, we will use the short notation
   $s_{\mu}$ to replace $u(\tilde{\pi}\times\lz)(s_{\mu})u^*$ in
   case it acts on $L^2(\Gamma, {\cal H})$, and then we have that
   for $\eta\in L^2(\Gamma, {\cal H})$,
   $$(s_{\mu}\eta)(\sigma)=s_{\mu}(\eta(\sigma-\oz_{\mu})),\;\ \
   \mbox{and}\;\ \ fs_{\mu}=s_{\mu}\sigma_{\om}(f)\; (\mbox{on}\; L^2(\Gamma, {\cal H})),$$
   where $f\in L^{\infty}(\Gamma)$,
   $\sigma_{\gz_0}(f)(\gz)=f(\gz +\gz_0)$, $\forall \gz, \gz_0 \in
   \Gamma$, and $f$ acts on $L^2(\Gamma, {\cal H})$ by the
   point-wise multiplication.
 Therefore for
  $f\in L^{\infty}(\Gamma)$, $f$ commutes with $s_{\mu}s_{\mu}^*$ and $p_v§$ for any $\mu\in E^*$ and $v\in E^0$. Moreover since
   $\widehat{C^*(G)}=C_0(\Gamma)$, we have that
   $$C^*(E)\times_{\az^{\oz}}G
   =\overline{span}\{s_{\mu}fs_{\nu}^*: \mu, \nu\in E^*, r(\mu)=r(\nu)\in E^0, f\in C_0(\Gamma)\}
   \subseteq {\bf B}(L^2(\Gamma, {\cal H}))$$

     For any subset $S$ of $\Gamma$, denote $\chi_S$ to be the
     characteristic function on $\Gamma$ of $S$. Let
     $\{U_i\}_{i\in\mathbf{I}}$ be an open base of $\Gamma$ such
     that for any $\gamma\in{\oz(E^*)}\subseteq \Gamma$,
     $j\in\mathbf{I}$, $\overline{U_j}$ compact and $U_j-\gamma\in\{U_i:\ i\in\mathbf{I}\}$, and
     let $\Lambda$ be the directed set consisting of all finite and not empty
     subsets of $\mathbf{I}$ with the inclusion order. Let
     $D_0(\Gamma)$ be the $C^*$-subalgebra of $L^{\infty}(\Gamma)$
     which is generated by all the characteristic functions
     $\chi_{U_i}(i\in\mathbf{I})$, and for any $\lz\in\Lambda$,
     let $D_{\lz}(\Gz)$ be the $C^*$-subalgebra of
     $L^{\infty}(\Gamma)$ which is generated by all the characteristic functions
     $\chi_{U_i}(i\in\lz)$. Let ${\cal F}(E)$ be the $C^*$-subalgebra of ${\bf B}(L^2(\Gamma, {\cal H}))$ generated
by the set $\{s_{\mu}fs_{\nu}^*: \mu, \nu\in E^*, f\in
D_0(\Gamma)\}$, and moreover let
 ${\cal F}_{\lz}(E)$ be the $C^*$-subalgebra of ${\bf B}(L^2(\Gamma, {\cal H}))$ generated
by the set $\{s_{\mu}fs_{\nu}^*: \mu, \nu\in E^*, f\in
D_{\lz}(\Gamma)\}$ which is equal to the set $\{s_{\mu}fs_{\nu}^*:
\mu, \nu\in E^*, r(\mu)=r(\nu)\in E^0, f\in D_{\lz}(\Gamma)\}$.
 Then we have the following relations easily.\\

{\sl{\bf Proposition 2.2.} With the notations as above, we have:\\
(1)
     $ C_0(\Gamma)\subseteq
     D_0(\Gamma)\;$  and  $\;C^*(E) \times_{\az^{\oz}}G\subseteq
     {\cal F}(E)$. If moreover $\Gamma$ is discrete (equivalently
     $G$ is compact), $\mathbf{I}=\Gamma$, and $U_i=\{i\}$ for any
     $i\in\mathbf{I}=\Gamma$, then
     $\;C^*(E) \times_{\az^{\oz}}G=
     {\cal F}(E)$.\\
(2)   $D_0(\Gamma)$ is the inductive limit of
    $D_{\lz}(\Gamma)$, and $D_0(\Gamma)$ is invariant under the actions of
    $\sigma_{\gamma}\ (\gamma\in{\oz(E^*)})$.\\
(3) ${\cal F}(E)$ is the inductive limit of ${\cal F}_{\lz}(E)$
with the coherent family of morphisms $\phi_{\lz_1,\lz_2}:
A_{\lz_1}\to A_{\lz_2}$ being inclusion map for ${\lz_1}\subseteq
{\lz_2}$.
 }\\

Since $D_0(\Gamma)$ is invariant under the actions of
$\sigma_{\gamma}\ (\gamma\in{\oz(E^*)})$, it is easy to see that
$$
{\cal F}(E)=\overline{span}\{s_{\mu}fs_{\nu}^*: \mu, \nu\in E^*,
r(\mu)=r(\nu)\in E^0,f\in D_0(\Gamma)\}\subseteq {\bf
B}(L^2(\Gamma,{\cal H})).
$$
 Since generally $C^*(E)$ may not be unital, it should be noted
 that $ C_0(\Gamma)$ is possibly not a subalgebra of ${\cal F}(E)$.

 From the Proposition 2.2, to embed $C^*(E)\times_{\az^{\oz}}G$ to an $AF$ algebra,
it is enough that ${\cal F}(E)$ is an $AF$ algebra. So in the rest
of this paper, we will mainly discuss ${\cal F}_{\lz}(E)$ to find
when ${\cal F}(E)$ is an $AF$ algebra. Therefore, for convenience,
we arbitrarily chose a $\lz\in\Lambda$, and fix this $\lz$ through
the rest of the paper except Lemma 3.7 and Theorem 3.8.

Since $D_{\lz}(\Gamma)$ is of finite dimension and abelian, with
the $\lz\in\Lambda$ fixed and mentioned above, there are mutually
orthogonal minimal projections $p_1, p_2, \cdots, p_M$ in
$D_{\lz}(\Gamma)$ such that $D_{\lz}(\Gamma)$ consists of all
their linear combinations. Let $p$ be the unit of
$D_{\lz}(\Gamma)$, then $p$ is the characteristic function of
$U=\cup_{i\in\lz}U_i$, and is the sum of all $p_i$. Then it is
easy to see that:
$$
{\cal F}_{\lz}(E) = \overline{alg.-span} \{s_{\mu}p_i s_{\nu}^*:
\mu, \nu\in E^*,
     r(\mu)=r(\nu)\in E^0, i=1,2,\cdots, M\}\subseteq {\bf B}(L^2(\Gamma, {\cal
     H})),
$$
 where for a set $X$ in a topological algebra $A$, $\overline{alg.-span}X$ is the
 closed subalgebra of $A$ generated by $X$, i.e. the smallest closed subalgebra of $A$
 which contains $X$. For convenience, we let $1=\chi_{\Gamma}$, which is
 the identity operator
on $ L^2(\Gamma, {\cal H})$, and let $p_0=1- p$.\\

Let ${\cal A}_{\lz}(E))$ be the (not closed) algebra generated by
$\{1\}$, $D_{\lz}(\Gamma)$ and $\{s_{\mu}fs_{\mu}^*:\ \mu\in E^*,
f\in \{1\}\cup D_{\lz}(\Gamma) \}$, then ${\cal A}_{\lz}(E))$ is a
$*$-subalgebra of $M({\cal F}(E)))$. For any $v\in E^0$ and $k\ge
1$, we define a map $\rho_{vk}$ on ${\cal A}_{\lz}(E)$ by
$$
\rho_{vk}(x)= \sum_{s(\mu)=v,|\mu|=k}s_{\mu}xs_{\mu}^* \ \
(\forall x\in{\cal A}_{\lz}(E)).
 $$
We note that if $v$ is a sink, then $\rho_{vk}=0$. Since $E$ is
row-finite,
the right hand of the equation above is a finite sum.\\

{\sl{\bf Lemma 2.3.} (1)  $\rho_{vk}$ is a a $*$-endmorphism on
${\cal A}_{\lz}(E)$, and so $\rho_{vk}(p)$ is a projection.

 (2) $\{s_{\mu}ps_{\mu}^*\}_{\mu\in E^*}$ are projections which commute with
each other, and so the projections $\{\rho_{vk}(p)\}_{v\in E^0,
k\ge 1}$ commute with each other too.

(3) $\rho_{vk}(p)p_u=p_u\rho_{vk}(p)=\delta_{u,v}\rho_{vk}(p)$ for
any $u,v\in E^0$ and $k\ge 1$. $s_{\mu}ps_{\mu}^*$ and
$\rho_{vk}(p)$ commute with $f$ for any $v\in E^0$, $k\ge 1$,
$\mu\in E^*$ and $f\in D_0(\Gamma)$.

(4) $\rho_{vk}(1-\rho_{ul}(p))=p_v$ for any $u,v\in E^0$, $k,l \ge
1$, if there is no path $\mu$ in $E$ from vertex $v$ to vertex $u$
with $|\mu|=k$.

(5) Let $v\in E^0$, $k\ge 1$, and $\Sigma$ be a subset of finite
set $\{\mu\in E^*: s(\mu)=v\mbox{ and }|\mu|=k\}$, then
$p_v-\sum\limits_{\mu\in\Sigma}s_{\mu}ps_{\mu}^*
=\prod\limits_{\mu\in\Sigma}(p_v-s_{\mu}ps_{\mu}^*)$, and
$1-\sum\limits_{\mu\in\Sigma}s_{\mu}ps_{\mu}^*
=\prod\limits_{\mu\in\Sigma}(1-s_{\mu}ps_{\mu}^*)$}.\\

{\bf Proof:} (1) Let $\mu\in E^*$, $ f\in D_{\lz}(\Gamma)$, $v\in
E^0$, then $(s_{\mu}fs_{\mu}^*)p_v
=\delta_{s(\mu),v}(s_{\mu}fs_{\mu}^*)=p_v (s_{\mu}fs_{\mu}^*)$.
Therefore for any $v\in E^0$, $p_v$ commutes with ${\cal
A}_{\lz}(E)$.

Let $x,y\in {\cal A}_{\lz}(E)$, then
$$
\begin{array}{l}
\rho_{vk}(x)\rho_{vk}(y)=(\sum\limits_{s(\mu)=v,|\mu|=k}s_{\mu}xs_{\mu}^*)
(\sum\limits_{s(\mu)=v,|\mu|=k}s_{\mu}ys_{\mu}^*)\\
\\
=\sum\limits_{s(\mu)=v,|\mu|=k}s_{\mu}xp_{r(\mu)}y s_{\mu}^*=
\sum\limits_{s(\mu)=v,|\mu|=k}s_{\mu}p_{r(\mu)}xy s_{\mu}^*\\
\\
=\sum\limits_{s(\mu)=v,|\mu|=k}s_{\mu}xys_{\mu}^*= \rho_{vk}(xy)
\end{array}
$$
Clearly $\rho_{vk}(x^*)=\rho_{vk}(x)^*$, i. e. $\rho_{vk}$ is a a
$*$-endmorphism on ${\cal A}_{\lz}(E)$.

(2) Since $p$ commutes with $p_v\ (\forall v\in E^0)$, for any
$\mu\in E^*$,
 $(s_{\mu}ps_{\mu}^*)^2=s_{\mu}ps_{\mu}^*$, i.e. $s_{\mu}ps_{\mu}^*$ is a projection.
 Therefore it is enough to prove that
$\{s_{\mu}ps_{\mu}^*\}_{\mu\in E^*}$ commute with each other. Let
$\mu,\nu\in E^*$, without loss of generality, we assume that
$|\nu|\ge |\mu|\ge 0$. If $\nu=\mu$, it is clear; If
$\nu=\mu\mu_1$ for some $\mu_1\in E^*\backslash E^0$, then
$$
\begin{array}{l}
(s_{\mu}ps_{\mu}^*)(s_{\nu}ps_{\nu}^*)=
s_{\mu}ps_{\mu_1}ps_{\mu\mu_1}^*\\
=s_{\mu\mu_1}ps_{\mu\mu_1}^*\sigma_{-\oz_{\mu}}(p)=
s_{\nu}ps_{\nu}^*s_{\mu}s_{\mu}^*\sigma_{-\oz_{\mu}}(p)\\
=(s_{\nu}ps_{\nu}^*)(s_{\mu}ps_{\mu}^*);
\end{array}
$$
In other case,
$(s_{\mu}ps_{\mu}^*)(s_{\nu}ps_{\nu}^*)=0=(s_{\nu}ps_{\nu}^*)(s_{\mu}ps_{\mu}^*)$.
Therefore $\{s_{\mu}ps_{\mu}^*\}_{\mu\in E^*}$ commute with each
other.

(3) is clear from direct computation.

(4) is from the fact that $s_\mu s_\nu=0$ for any $\mu,\nu\in E^*$
with $s(\mu)=v$, $|\mu|=k$ and $s(\nu)=u$, since $r(\mu)\not=u$
under the assumption.

(5) is from that $(s_{\mu}ps_{\mu}^*)(s_{\nu}ps_{\nu}^*)=0$ for
any $\mu,\nu\in\Sigma$ with $\mu\not=\nu$, and
$p_v(s_{\mu}ps_{\mu}^*)=s_{\mu}ps_{\mu}^*$  for
any $\mu\in E^*$ with $s(\mu)=v$.\\

Certainly it is impossible to embed all
$C^*(E)\times_{\az^{\oz}}G$ defined above to an $AF$ algebra. For
example, if $C^*(E)$ has an infinite projection and $G$ is
discrete,
  $C^*(E)\times_{\az^{\oz}}G$ is not a subalgebra of an $AF$ algebra.
  The following proposition gives a necessary condition for
  embedding  $C^*(E)\times_{\az^{\oz}}G$ to an $AF$ algebra in another special case.\\

{\sl{\bf Proposition 2.4.} Let $E$ be a row-finite directed graph,
$G$ be a locally compact abelian group with dual $\hat{G}=\Gamma$,
$\oz$ be a labelling map from $E^*$ to $\Gamma$ such that
$C^*(E)\times_{\az^{\oz}}G$ has no infinite projection ( specially
if $C^*(E)\times_{\az^{\oz}}G$ can be embedded to an $AF$
algebra). Then $\oz_{\gamma}+O\not= O$ for any loop $\gamma$ in
$E$ which has an exit, and any compact open subset $O$ of
$\Gamma$. In particular if  $\Gamma$ has a compact open subset (
specially if $G$ is compact ), then
 $\oz_{\gamma}\not=1_{\Gamma}$ for any loop $\gamma$ in $E$ which has an
exit.\\
}

{\bf Proof:} Since $O$ be a compact open subset of $\Gamma$, then
$\chi_O$ is a projection in $C_0(\Gamma)$. If there is a loop
$\gamma$ in $E$ with an exit such that $\oz_{\gamma}+O= O$, let
$v$ be a vertex in $\gamma$ at which there is an exit. Viewing
$\gamma$ as a path with $s(\gamma)=r(\gamma)=v$, we have
$p_v\chi_O$, $s_{\gamma}\chi_O s_{\gamma}^* \in
C^*(E)\times_{\az^{\oz}}G$. Moreover $s_{\gamma}\chi_O
s_{\gamma}^* =\sigma_{-\oz_{\gamma}}(\chi_O)s_{\gamma}
s_{\gamma}^*=\chi_{(\oz_{\gamma}+O)} s_{\gamma}
s_{\gamma}^*=\chi_{O} s_{\gamma} s_{\gamma}^*$, and $(\chi_O
s_{\gamma}^*)(s_{\gamma}\chi_O)=\chi_O p_v$, i.e. $\chi_O
s_{\gamma} s_{\gamma}^*\sim \chi_O p_v$. Since clearly $\chi_O
s_{\gamma} s_{\gamma}^*< \chi_O p_v$ for there is an exit of
$\gamma$ at $v$, i.e. $ \chi_O p_v$ is a infinite projection in
$C^*(E)\times_{\az^{\oz}}G$, which contradicts that
$C^*(E)\times_{\az^{\oz}}G$ has no infinite projection.
\\

In our $AF$ embedding frame, we will in fact embed the bigger
$C^*$-algebra ${\cal F}(E)$, which contains
$C^*(E)\times_{\az^{\oz}}G$ and ${\cal F}_{\lz}(E)$ by Proposition
2.2(1) and definitions, into an $AF$ algebra. By the similar proof
as that of Proposition 2.4, without the assumption that $\Gamma$
has a compact open subset, we still have
$\oz_{\gamma}\not=1_{\Gamma}$ for any loop $\gamma$ in $E$ which
has an exit, if ${\cal F}_{\lz}(E)$ can be embedded into an $AF$
algebra (Since then there is an infinite projection $ \chi_{O}p_v$
in ${\cal F}_{\lz}(E)$ for $v\in\gamma$ and $O=U_i$ with
$i\in\lz$, if $\oz_{\gamma}=1_{\Gamma}$ for some loop $\gamma$ in
$E$ ).

Now to embed $C^*(E)\times_{\az^{\oz}}G$ to an $AF$ algebra, first
we
want to define the labelling map, which defines the action, to be almost proper.\\

{\sl{\bf Definition 2.5.} Let $\oz$ be a labelling map from $E^*$
to $\Gamma$. With the discrete topology, $E^*$ is a topological
space. We call $\oz$ to be almost proper if $\oz|_{E^*\backslash
E^0}$ is proper, i.e. for any compact
subset $A$ of $\Gamma$, $\oz^{-1}(A)\backslash E^0$ is a finite set.}\\

It is clear that if $E^*\backslash E^0$ is an infinite set and
$\oz$ is almost proper, then $\Gamma$ is not a compact set, which
is equivalent to that $G$ is not discrete; and if $E$ has a loop
$\gamma$ and $w$ is almost proper, then
$w_{\gamma}\not=1_{\Gamma}$, for otherwise
$w^{-1}(\{1_{\Gamma}\})\backslash E^0$ is an infinite set. The
following proposition says that this definition can be viewed as a
reasonable generalization of the corresponding definition given in
 {\cite{[Ka]}}.\\

{\sl{\bf Proposition 2.6.} Let $E=(E^0,E^1,r,s)$ be a finite
directed graph such that there is a $v\in E^0$ with $r(e)=s(e)=v$
for any $e\in E^1$, $\oz$ be a labelling map from $E^*$ to
$\Gamma$. Then $-\oz_e\notin \overline{\oz(E^*)}$ if and only if
$\oz$ is almost proper. } \\

{\bf Proof:} If $-\oz_e\notin \overline{\oz(E^*)}$ for each $e\in
E^1$, and $\oz$ is not almost proper, then there is a compact
subset $A$ of $\Gamma$ such that $\oz^{-1}(A)\backslash E^0$ is
not finite. Since $E^1$ is finite and $A$ is compact, it is easy
to chose a sequence $\mu_1$, $\mu_2$, $\cdots$, $\mu_i$, $\cdots$
in $\oz^{-1}(A)\backslash E^0$ such that $\oz_{\mu_i}$ is
convergent in $\Gamma$ and $\lim_{i\to\infty}|\mu_i|=+\infty$.
Moreover by the finiteness of $E^1$ again, there are $f\in E^1$
and a subsequence $\{\mu_{k_i}\}$ of $\{\mu_i\}$ such that
$\{l_i(e)\}_{i\ge 1}$ is increasing for any $e\in E^1$, and
$\lim_{i\to\infty}l_i(f)=+\infty$, where $l_i(e)$ is the appearing
times of $e$ in $\{\mu_{k_i}\}$. Without loss of generality, we
may assume $\{l_i(f)\}_{i\ge 1}$ is monotonous strictly, and then
$-\oz_f=\lim_{i\to\infty}(\oz_{\mu_{k_{i+1}}}-\oz_{\mu_{k_{i}}}-\oz_{f})\in
\overline{\oz(E^*)}$, which is a contradiction.

For the converse, if $\oz$ is almost proper, and $-\oz_f\in
\overline{\oz(E^*)}$ for some $f\in E^1$, then
$-\oz_f=\lim_{i\to\infty}\oz_{\mu_i}$ for some $\{\mu_i\}\subseteq
E^*$. Since $r(e)=s(e)=v$ for any $e\in E^1$, let $nf$ be the path
with edge $f$ repeating $n$ times, then
$-\oz_{nf}=n(-\oz_f)=\lim_{i\to\infty}n\oz_{\mu_i}
=\lim_{i\to\infty}\oz_{n\mu_i}\in \overline{\oz(E^*)}$. Let ${\cal
V}$ be the neighborhood base of $1_{\Gamma}$ in $\Gamma$, and let
$\Oz$ be the directed set $\{t=(n,O): n\in{\bf N}, O\in{\cal V}\}$
with the order: $t_1=(n_1,O_1)\leq t_2=(n_2,O_2)$ iff $n_1\leq
n_2$ and $O_2\subseteq O_1$. For any $t=(n,O)$, since
$-\oz_{nf}\in \overline{\oz(E^*)}$, we may chose a $\mu\in E^*$
such that
 $\oz_{nf}+\oz_{\mu}\in O$. Let $\mu_{t}=(nf)\mu\in E^*$,
 then $|\mu_{t}|\geq n$ and $\oz_{\mu_{t}}\in O$, i.e.
 $\lim_{t}|\mu_{t}|=+\infty$ and $\lim_{t}\oz_{\mu_{t}}=1_{\Gamma}$. This
contradicts that $\oz$ is almost proper.\\

{\bf Example 2.7.} (1) If $E$ is a finite directed graph without
loop (specially if $E$ is a finite rooted directed tree), then
$E^*$ is finite, and so any labelling map is almost proper.

(2) Let $G=\Gamma={\bf R}$, $E=(E^0, E^1, r, s)$ be a directed
graph with $E^1=\{e_1, e_2, \cdots \}$, $\oz: E^*\to {\bf R}$ be a
labelling map with $\oz_{e_i}\oz_{e_j}>0\ (\forall i, j)$. If
$E^1$ is finite or $\lim_{n\to +\infty}\oz_{e_n}=\infty$, then $\oz$ is almost proper.\\

{\sl{\bf Lemma 2.8.} Let $p$ be the unit of $D_{\lz}(\Gamma)$,
$\oz$ be a labelling map from $E^*$ to $\Gamma$. If $\oz$ is
almost proper, then there is a $N\in{\bf N}\cup\{0\}$, and
$F=\{\mu_1, \mu_2,\cdots,\mu_N\}\subseteq E^*\backslash E^0$ such
that $ps_{\mu}p=0$ for any $ \mu\in E^*\backslash (E^0\cup F)$. }\\

{\bf Proof:} Let $F=\oz^{-1}(\overline{U-U})\backslash E^0$ with
$U=\cup_{i\in\lz}U_i$, since $\overline{U-U}$ is compact, and
$\oz$ is almost proper, $F$ is finite. There are $N\in{\bf
N}\cup\{0\}$ and $\mu_1, \mu_2,\cdots,\mu_N\in E^*\backslash E^0$
such that $F=\{\mu_1, \mu_2,\cdots,\mu_N\}$. For any $\mu\in E^*$,
if $ps_{\mu}p\not=0$, then $0\not=s_{\mu}^* p
s_{\mu}p=\sigma_{\oz_{\mu}}(p)ps_{\mu}^*
s_{\mu}=\sigma_{\oz_{\mu}}(p)pp_{r(\mu)}$, and so there is a $x\in
U$ such that $x+\oz_{\mu}\in U$. Therefore $\oz_{\mu}\in
U-U\subseteq \overline{U-U}$, i.e.
$\mu\in F\cup E^0$, and this completes the proof.\\

\section*{3. $AF$-embedding}

\hspace{0.6cm} Let $T=(T^0, T^1, r, s)$ be a row-finite directed
graph whose each weak connected component is a rooted (possibly
infinite) directed tree, and let $E=(E^0, E^1, r, s)$ be
 the row-finite directed graph constructed from $T$ by attaching
 $n_v$ ($0\le n_v <+\infty$) 1-loops to each vertex $v$ in $T$. It is clear
that the directed graph with one vertex and $n$ 1-loops, whose
$C^*$-algebra is the Cuntz algebra, is a special one of these $E$.
In this section we always assume $E$ to be of this form and the
labelling map $\oz: E^*\to \Gamma$ to be almost proper.

From \cite{[BPRS]} or \cite{[F3]}, it is easy to see that the
$C^*$-algebras $C^*(E)$ of the graphs under consideration can be
simple or non-simple, and also can be purely infinite or finite
(see some simple examples below). Moreover among them there are
some interesting $C^*$-algebras. The reason only to consider this
class of graphs here is that the pathes in a general graph is too
complex to handle in our embedding construction.
\\

{\bf Example 3.1.} Let $T=(T^0, T^1, r, s)$ be a (possibly
infinite) rooted directed tree with $|T^0|=|T^1|+1$, where
$T^0=\{v_i:1\le i<|T^0|+1 \}$ and $T^1=\{e_i:1\le i<|T^1|+1 \}$
are two countable (possibly infinite) sets with $s(e_i)=v_i$,
$r(e_i)=v_{i+1}$ for all $1\le i<|T^1|+1 $. Let $t=|T^1|$.

(1) If $1\le t\le +\infty$, $E$=$T$, then $C^*(E)$ is a simple
$C^*$-algebra with each $p_{v_n}\ (n\ge 1)$ finite;

(2)  If $t< +\infty$, $E$ is constructed from $T$ by attaching
finite and more than one 1-loops to $v_{t+1}$, then $C^*(E)$ is a
purely infinite simple $C^*$-algebra;

(3)  Let $\{l_i\}_{i\in{\bf N}}$ be a strictly increasing positive
integer sequence with $\lim_{i\to +\infty}l_i =+\infty$. If
$t=+\infty$ and $E$ is constructed from $T$ by attaching finite
1-loops to each $v_{l_i}\ (i\ge 1)$, then $C^*(E)$ is a purely
infinite non-simple $C^*$-algebra;

(4)  If $1\le t\le +\infty$, $E$ is constructed from $T$ by
attaching finite 1-loops to some $v_{n}\ (1\le n<t)$, then
$C^*(E)$ is a non-simple $C^*$-algebra with
each $p_{v_i}\ (i\le n)$ infinite and each $p_{v_j}\ (j>n)$ finite.\\

For any $v\in E^0$, $k\ge 1$, $0\le l\le k$, let
$$
E^{kl}_v=\{\mu=e_1 e_2\cdots e_k\in E^k: s(\mu)=v\mbox{ and }e_1,
e_2, \cdots, e_l \mbox{ are  1-loops, but }e_{l+1} \mbox{ is not a
1-loop}\},$$
$$
\rho^l_{vk}(x) =\sum\limits_{\mu\in E^{kl}_v}s_{\mu}xs_{\mu}^*\ \
(\forall x\in{\cal A}_{\lz}(E)).
$$
By the similar proof as that of Lemma 2.3(1), $\rho^l_{vk}$ is
still a homomorphism. Let $E^{k}_v=\{\mu\in E^k: s(\mu)=v\}$, then
$E^{k}_v=\cup_{l=0}^k E^{kl}_v$, and therefore
$\rho_{vk}=\sum_{l=0}^k \rho^l_{vk}$.

Let  $F=\{\mu_1, \mu_2,\cdots,\mu_N\}$ be defined in Lemma 2.8,
and let $W=\{s(\mu_1), s(\mu_2),\cdots,s(\mu_N)\}$, then
$(W,\le_T)$ is a partially ordered set. Let $V$ be the subset of
$W$ consisting of all the maximal elements in $(W,\le_T)$, and let
$$
m=\max\limits_{1\le i\le N}|\mu_i|+\max\{|\mu|: \mu \mbox{ is a
path in }T\mbox{ from a vertex in }V \mbox{ to a vertex in } W\}.
$$
For $v\in V$, $j\ge 0$, let $V(v,j)$ be all the vertices to which
there is a (unique) path $\mu$ in $T$ from $v$ with $|\mu|=j$, and
let
$$
E^*_F=\bigcup_{v\in V}\bigcup_{j=0}^{m-1}\bigcup_{u\in V(v,j)}
\{\mu: s(\mu)=u, 1\le|\mu|\le m-j\}.
$$
We note that if $F=\emptyset$, then $E^*_F=\emptyset$. It is clear
that for any $\mu_i\in F$, there are $v\in V$ and $j\in{\bf N}$
with $0\le j\le m-\max_{1\le i\le N}|\mu_i|$ such that
$s(\mu_i)\in V(v,j)$, and so $\mu_i \in E^*_F$. Since each weak
connected component of $T$ is a rooted directed tree, for any
$u,v\in V$ with $u\not= v$, there is no path in $E$ from one
vertex in $V(u,i)$ to another vertex in $V(v,j)$ for any $i,j\ge
0$, and also  there is no path in $E$ from one vertex in $V(v,k)$
to another vertex in $V(v,l)$ for any $k>l$. Now we let
$$
q=\prod\limits_{\mu\in E^*_F}(1-s_{\mu}ps_{\mu}^*)p,
$$
where $p$ is the unit of $D_{\lz}(\Gamma)$ for the given $\lz\in
\Lambda$ as mentioned in part 2. We note that if $F=\emptyset$,
then $q=p$. By Lemma 2.3 (2) (3), the definition above is
well-defined and $q\le p$. By Lemma 2.3 (5),
$1-\rho_{vk}(p)=\prod\limits_{s(\mu)=v,
|\mu|=k}(1-s_{\mu}ps_{\mu}^*)$, and so
$$
\begin{array}{l}
q=\prod\limits_{v\in V}\prod\limits_{j=0}^{m-1}\prod\limits_{u\in
V(v,j)} \prod\limits_{k=1}^{m-j}
(1-\rho_{uk}(p))p\\
=\prod\limits_{v\in V}(\prod\limits_{k=1}^{m}(1-\rho_{vk}(p))
\prod\limits_{u\in V(v,1)}\prod\limits_{k=1}^{m-1}(1-\rho_{uk}(p))
\cdots \prod\limits_{u\in V(v,m-1)}(1-\rho_{u1}(p)))p,
\end{array}
$$
which is a projection in ${\cal A}_{\lz}(E)\subseteq M({\cal F}(E))\subseteq {\bf B}(L^2(\Gamma,{\cal H}))$ by Lemma 2.3 (1).\\

{\sl{\bf Lemma 3.2.} With the notations as above, we have:

(1) $q$ and $\rho_{uk}(q)$ commute with $f$ and $p_v$ for any
$u,v\in E^0$, $k\ge 1$ and $f\in D_0(\Gamma)$.

(2) For any $\mu,\nu\in E^*$, $(s_{\mu}qs_{\mu}^*
)(s_{\nu}qs_{\nu}^*)= \delta_{\mu,\nu}s_{\mu}qs_{\mu}^*$. If
moreover $\mu\notin E^0$, $qs_{\mu}q=0$.
}\\

{\bf Proof:} (1) is clear from Lemma 2.3 and direct computation.

(2) is also clear if $\mu,\nu\in E^0$. Since
$(s_{\mu}qs_{\mu}^*)^2=s_{\mu}qs_{\mu}^*$ by (1), it is enough to
prove that for any $\mu\in E^*\backslash E^0$, $qs_{\mu}q=0$,
which is equivalent to
 $qs_{\mu}qs_{\mu}^* q=0$. If $\mu\notin F$, by Lemma 2.8,
$qs_{\mu}qs_{\mu}^*q=q(ps_{\mu}p)qs_{\mu}^*q=0$. If $\mu\in F$, by
the discussion above $\mu\in E^*_F$. Therefore
$0\le qs_{\mu}qs_{\mu}^*q\le q(s_{\mu}ps_{\mu}^*)q =0$, i.e. $qs_{\mu}qs_{\mu}^* q=0$. \\

{\sl{\bf Lemma 3.3.} With the notations as above, let $v\in E^0$,
then there is a finite subset $E^*(v)$ of $E^*$ such that
$$
p_v p=\sum\limits_{\mu\in E^*(v)}s_{\mu}qs_{\mu}^* p.
$$
Moreover if $v\notin\bigcup\limits_{v\in V, 0\le j\le m-1}V(v,j)$,
$E^*(v)=\{v\}$; and if $v\in\bigcup\limits_{v\in V, 0\le j\le
m-1}V(v,j)$, $E^*(v)$ is a subset of  finite set
$E^*_F\bigcup(\bigcup\limits_{v\in V, 0\le j\le m-1}V(v,j))$.
}\\

{\bf Proof:}  By definition, $q=\prod\limits_{v\in
V}\prod\limits_{j=0}^{m-1}\prod\limits_{u\in
V(v,j)}\prod\limits_{k=1}^{m-j} (1-\rho_{uk}(p))p$.

Then if $v\notin\bigcup_{v\in V, 0\le j\le m-1}V(v,j)$, it is
clear from Lemma 2.3 (3) that $p_v qp_v p=p_v p$.

If now $v\in V$, by Lemma 2.3 (3),
%$$
\begin{equation}
p_v qp_v p=p_v \prod\limits_{k=1}^m(1-\rho_{vk}(p))p_v p
=\prod\limits_{k=1}^m(p_v -\rho_{vk}(p))p.
\end{equation}
%$$
Let $1\le k\le m$, since there is no path from $v$ to any vertex
in $V(w,j)$ for any $0\le j\le m-1$ and $w\in V$ with $w\not= v$,
by Lemma 2.3 (4),
$$
\begin{array}{l}
\rho_{vk}(q)p=\rho_{vk}(\prod\limits_{j=0}^{m-1}\prod\limits_{u\in
V(v,j)} \prod\limits_{s=1}^{m-j}
(1-\rho_{us}(p))p)p\\
\\
=\prod\limits_{j=0}^{m-1}\prod\limits_{u\in V(v,j)}
\prod\limits_{s=1}^{m-j}
(p_v-\rho_{vk}(\rho_{us}(p)))\rho_{vk}(p)p\ \  (\mbox{ for
}\rho_{vk}
\mbox{ is a homomorphism})\\
\\
=\prod\limits_{j=0}^{m-1}
\prod\limits_{s=1}^{m-j}\prod\limits_{u\in V(v,j)}
(p_v-\rho_{vk}(\rho_{us}(p)))\rho_{vk}(p)p
\end{array}
$$
Since $\mu\notin F(\subseteq E^*_F)$ for any $\mu\in E^*$ with
$s(\mu)\in V(v,j)$ ( $\forall v\in V$, $0\le j\le m-1$ ) and
$|\mu|>m-j$, we have
$$
\begin{array}{l}
\rho_{vk}(q)p=\prod\limits_{j=0}^{m-1}
\prod\limits_{s=1}^{m}\prod\limits_{u\in V(v,j)}
(p_v-\rho_{vk}(\rho_{us}(p)))\rho_{vk}(p)p \\
(\mbox{ since }\rho_{vk}(\rho_{us}(p))p=0\mbox{ for }s>m-j\mbox{
and }
u\in V(v,j)\mbox{ by Lemma 2.8 })\\
\\
=\prod\limits_{s=1}^{m}\prod\limits_{j=0}^{m-1} \prod\limits_{u\in
V(v,j)} (p_v-\rho_{vk}(\rho_{us}(p)))\rho_{vk}(p)p
\end{array}
$$
$$
\begin{array}{l}
=\prod\limits_{s=1}^{m-k}\prod\limits_{j=0}^{m-1}
\prod\limits_{u\in V(v,j)}
(p_v-\rho_{vk}(\rho_{us}(p)))\rho_{vk}(p)p \\
(\mbox{ since }\rho_{vk}(\rho_{us}(p))p=0\mbox{ for }s>m-k\mbox{
and }
u\in V(v,j)\mbox{ by Lemma 2.8 })\\
\\
=\prod\limits_{s=1}^{m-k}\prod\limits_{j=0}^{m-1}
\prod\limits_{u\in V(v,j)}
(p_v-\sum\limits_{l=0}^k\rho_{vk}^l(\rho_{us}(p)))\rho_{vk}(p)p \\
\\
=\prod\limits_{s=1}^{m-k}(p_v-\sum\limits_{j=0}^{m-1}
\sum\limits_{u\in V(v,j)}
\sum\limits_{l=0}^k\rho_{vk}^l(\rho_{us}(p)))\rho_{vk}(p)p \ \
(\mbox{ by Lemma 2.3 (5) })\\
\\
=\prod\limits_{s=1}^{m-k}(p_v-\sum\limits_{l=0}^k
\sum\limits_{j=0}^{m-1} \sum\limits_{u\in V(v,j)}
\rho_{vk}^l(\rho_{us}(p)))\rho_{vk}(p)p.
\end{array}
$$
For $1\le s\le m-k$, since
$$
\sum\limits_{j=0}^{m-1} \sum\limits_{u\in V(v,j)}
\rho_{vk}^k(\rho_{us}(p))=\rho_{vk}^k(\rho_{vs}(p))
=\sum\limits_{l=k}^{k+s}\rho_{v(k+s)}^l(p),
$$
and for $0\le l\le k-1$,
$$
\sum\limits_{j=0}^{m-1} \sum\limits_{u\in V(v,j)}
\rho_{vk}^l(\rho_{us}(p))=\rho_{v(k+s)}^l(p),
$$
we have
$$
\begin{array}{l}
\rho_{vk}(q)p
=\prod\limits_{s=1}^{m-k}(p_v-\sum\limits_{l=0}^{k+s}
\rho_{v(k+s)}^l(p))
\rho_{vk}(p)p\\
\\
=\prod\limits_{s=1}^{m-k}(p_v-\rho_{v(k+s)}(p)) \rho_{vk}(p)p.
\end{array}
$$
And so
\begin{equation}
\rho_{vk}(q)p=\prod\limits_{s=k+1}^{m}(p_v-\rho_{vs}(p))
\rho_{vk}(p)p.
\end{equation}

Therefore $p_v p =p_v q p_v p+\sum\limits_{k=1}^m \rho_{vk}(q)p$
by (1), (2) and direct computation.\\

If  $v\in V(u,i)$ for some $u\in V$ , $1\le i\le m-1$, the
discussion is similar only by replacing $m$ with $m-i$ and
replacing $V(v,j)$ ( $0\le j \le m-1$ ) with $\widetilde{V(v,j)}$
( $0\le j \le m-i-1$ ), where $\widetilde{V(v,j)}$ consists of all
vertices to which there is a unique
path $\mu$ in $T$ from $v$ with $|\mu|=j$.\\

Let $p_0=1-p$, and let $p_1,p_2, \cdots, p_M$ be the mutually
orthogonal minimal projections in $D_{\lz}$ with $\sum_{i=1}^M
p_i=1$, which are defined in part 2. We arbitrarily chose a vertex
$\tilde{v}\in E^0$, and let
$\widetilde{E^*_F}=E^*_F\cup\{\tilde{v}\}$. Let ${\bf J}$ be the
finite set of all maps from $\widetilde{E^*_F}$ to the set $\{0,
1, 2, \cdots, M\}$, and ${\bf K}_1=\{(v,\tau): v\in E^0, \tau\in
{\bf J}\}$. For $\forall (v,\tau)\in {\bf K}_1$, Let
$$
q_{(v,\tau)}=q\prod_{\mu\in
\widetilde{E^*_F}}\sigma_{w_{\mu}}(p_{\tau(\mu)})p_v\in{\bf
B}(L^2(\Gamma, {\cal H})).
$$
Clearly $q_{(v,\tau)}$ is a projection with $q_{(v,\tau)}\le q\le
p$.
 Let ${\bf K}=\{(v, \tau)\in {\bf K}_1: q_{(v,\tau)}\not=0\}$,
 which is clearly a countable (possibly infinite) set.\\

{\sl{\bf Lemma 3.4.} (1) $\{q_{(v,\tau)}\}_{(v,\tau)\in {\bf K}}$
are mutually orthogonal projections in ${\cal F}(E)$.

(2) For any fixed $v\in E^0$, $p_v q=\sum\limits_{\tau\in{\bf
J}}q_{(v,\tau)} =\sum\limits_{\tau:(v,\tau)\in{\bf
K}}q_{(v,\tau)}$ }

{\bf Proof:} (1) It is clear that $\{q_{(v,\tau)}\}_{(v,\tau)\in
{\bf K}}$ are mutually orthogonal. Since $D_0(\Gamma)$ is
invariant under the action of $\sigma_{\gamma}\
(\gamma\in{\oz(E^*)})$, $\prod\limits_{\mu\in
\widetilde{E^*_F}}\sigma_{w_{\mu}}(p_{\tau(\mu)})p_v p\in {\cal
F}(E)$. Therefore
 $q_{(v,\tau)}=q\prod\limits_{\mu\in \widetilde{E^*_F}}\sigma_{w_{\mu}}(p_{\tau(\mu)})p_v p
 \in {\cal F}(E)$ for $q\in {\cal A}_{\lz}(E)\subseteq M({\cal F}(E))$.\\

(2) $$
\begin{array}{l}
p_v q=p_v q \prod\limits_{\mu\in \widetilde{E^*_F}}\sigma_{w_{\mu}}(p_0+p_1+\cdots p_M)\\
\\
=p_v q\sum\limits_{\tau\in{\bf J}} \prod\limits_{\mu\in
\widetilde{E^*_F}}\sigma_{w_{\mu}}(p_{\tau(\mu)})
=\sum\limits_{\tau\in{\bf J}}q_{(v,\tau)}
\end{array}
$$
\\
{\sl{\bf Lemma 3.5.} The $C^*$-algebra ${\cal F}_{\lz}$ generated
by $\{s_{\mu}q_{(v,\tau)}s_{\nu}^*\}_{\mu,\nu\in E^*,
(v,\tau)\in{\bf K}}$ is isomorphic to $\bigoplus_{(v,\tau)\in{\bf
K}}{\cal K}_{(v,\tau)}$, and hence is an $AF$ subalgebra of ${\cal
F}(E)$, where ${\cal K}_{(v,\tau)}$ is a compact operator algebra
over a finite dimensional or a separable infinite dimensional
Hilbert space.
}\\

{\bf Proof:} First since $q_{(v,\tau)}\in {\cal F}(E)$,
$s_{\mu}q_{(v,\tau)}s_{\nu}^*=(s_{\mu}p)q_{(v,\tau)}(s_{\nu}p)^*\in
{\cal F}(E)$. Let $(v_1,\tau_1)$, $(v_2,\tau_2)\in{\bf K}$, by
Lemma 3.2, for any $\mu_1,\mu_2,\nu_1,\nu_2\in E^*$,
$$
\begin{array}{l}
(s_{\mu_1}q_{(v_1,\tau_1)}s_{\nu_1}^*)(s_{\mu_2}q_{(v_2,\tau_2)}s_{\nu_2}^*)
=s_{\mu_1}q_{(v_1,\tau_1)}p_{v_1}(q s_{\nu_1}^* s_{\mu_2}q)q_{(v_2,\tau_2)}s_{\nu_2}^*\\
=\delta_{\nu_1,\mu_2} \delta_{v_1,r(\nu_1)}
\delta_{(v_1,\tau_1),(v_2,\tau_2)}
s_{\mu_1}q_{(v_1,\tau_1)}s_{\nu_2}^*.
\end{array}
$$

Let $(v,\tau)\in{\bf K}$. If $r(\mu)$ or $r(\nu)$ is not $v$, then
$s_{\mu}q_{(v,\tau)}s_{\nu}^*=0$. Let $E^{*,v}=\{\mu\in E^*:\
r(\mu)=v\}$, by the computation above,
$\{s_{\mu}q_{(v,\tau)}s_{\nu}^*\}_{\mu,\nu\in E^{*,v}}$ is a
matrix unit, and generates the same $C^*$-algebra, notated by
${\cal F}_{(v,\tau)}$, as
$\{s_{\mu}q_{(v,\tau)}s_{\nu}^*\}_{\mu,\nu\in E^*}$. Moreover
${\cal F}_{(v,\tau)}$ is isomorphic to a compact operator algebra
${\cal K}_{(v,\tau)}$ over a finite dimensional or a separable
infinite
dimensional Hilbert space. This completes the proof.\\

{\sl{\bf Lemma 3.6.} The $C^*$-algebra ${\cal F}_{\lz}(E)$ is a
subalgebra of
${\cal F}_{\lz}$.}\\

{\bf Proof:} Since ${\cal F}_{\lz}(E) = \overline{alg.-span}
\{s_{\mu}p_i s_{\nu}^*: \mu, \nu\in E^*,
     r(\mu)=r(\nu)\in E^0, i=1,2,\cdots,M\}$, it is enough to prove for each $1\le i\le M$,
     $\mu,\nu\in E^*$ with $r(\mu)=r(\nu)$, $s_{\mu}p_i s_{\nu}^*\in
{\cal F}_{\lz}$. For any $v\in E^0$, $1\le i\le M$,
$$
\begin{array}{l}
p_v p_i=p_v pp_i=\sum\limits_{\gamma\in
E^*(v)}s_{\gamma}qs_{\gamma}^* p_i\
(\mbox{ by Lemma 3.3 })\\
\\
=\sum\limits_{\gamma\in E^*(v)}s_{\gamma}(p_{r(\gamma)}q)s_{\gamma}^* p_i\\
\\
=\sum\limits_{\gamma\in E^*(v)}s_{\gamma}(\sum\limits_{\tau\in{\bf
J}}q_{(r(\gamma),\tau)})s_{\gamma}^* p_i\
(\mbox{ by Lemma 3.4 (2) })\\
\\
=\sum\limits_{\gamma\in E^*(v)}\sum\limits_{\tau\in{\bf J}}
s_{\gamma}q_{(r(\gamma),\tau)}\sigma_{\oz_{\gamma}}(p_i)s_{\gamma}^* \\
\\
=\sum\limits_{\gamma\in E^*(v)}\sum\limits_{\tau\in{\bf J}}
s_{\gamma}q\prod\limits_{\mu\in E^*_F}\sigma_{w_{\mu}}(p_{\tau(\mu)})
\sigma_{\oz_{\gamma}}(p_i)p_{r(\gamma)}s_{\gamma}^* \\
\\
=\sum\limits_{\gamma\in E^*(v)\cap E_F^*}\sum\limits_{\tau\in{\bf
J}:\tau(\gamma)=i}s_{\gamma}q_{(r(\gamma),\tau)}s_{\gamma}^*
+\sum\limits_{\gamma\in E^*(v)\backslash
E_F^*}\sum\limits_{\tau\in{\bf
J}:\tau(\tilde{v})=i}s_{\gamma}q_{(r(\gamma),\tau)}s_{\gamma}^*\\
(\mbox{since }  E^*(v)\backslash E_F^*\subseteq E^0, \mbox{ and }
\omega_v=1_{\Gamma} \mbox{ for any }  v\in E^0) .
\end{array}
$$
Let $r(\mu)=r(\nu)=v$, then
$$
\begin{array}{l}
 s_{\mu}p_i s_{\nu}^*=s_{\mu}p_v p_i
s_{\nu}^*\\
\\=\sum\limits_{\gamma\in E^*(v)\cap E_F^*}
\sum\limits_{\tau\in{\bf
J}:\tau(\gamma)=i}s_{\mu\gamma}q_{(r(\gamma),\tau)}s_{\nu\gamma}^*
 +\sum\limits_{\gamma\in E^*(v)\backslash
E_F^*}\sum\limits_{\tau\in{\bf
J}:\tau(\tilde{v})=i}s_{\mu\gamma}q_{(r(\gamma),\tau)}s_{\nu\gamma}^*\in{\cal
F}_{\lz}.
\end{array}
$$
This completes the proof.\\

{\sl{\bf Lemma 3.7.} Let $\Omega$ be a directed set,
$\{A_{\lz}\}_{\lz\in\Omega}$ be $C^*$-algebras, and the separable
$C^*$-algebra
$A=\lim\limits_{\longrightarrow}(A_{\lz},\phi_{\lz\mu})$ be the
inductive limits with the coherent family of morphisms
$\phi_{\lz\mu}: A_{\lz}\to A_{\mu}(\lz<\mu)$
being the inclusion maps. If each $A_{\lz}$ is contained in an $AF$
subalgebra of $A$, then $A$ is an $AF$ algebra.}\\

{\bf Proof:} This is a direct consequence of the local
characterization of $AF$ algebra (for example see \cite{[Da]}
Theorem III.3.4). We note that $A_{\lz}$ can be viewed as a
subalgebra of $A$ for $\phi_{\lz\mu}$ are the inclusion maps.\\

{\sl{\bf Theorem 3.8.} Let $T$ be a row-finite directed graph
whose each weak connected component is a rooted (possibly
infinite) directed tree, $E$ be a row-finite directed graph
constructed by attaching $n_v$ ($0\le n_v <+\infty$) 1-loops to
each vertex $v$ in $T$, $G$ be a locally compact abelian group
with the dual $\Gamma$, $\oz: E^*\to \Gamma$ be an almost proper
labelling map, $(C^*(E),G,\az^{\oz})$ be the $C^*$-dynamical
system with $\az^{\oz}$ being the quasi-free action defined by
$\oz$. Then the crossed product $C^*(E)\times_{\az^{\oz}}G$ can be
embedded into an $AF$ algebra. If moreover $G$ is compact, then
the crossed product
$C^*(E)\times_{\az^{\oz}}G$ itself is an $AF$ algebra.}\\

{\bf Proof:} In Lemma 3.7, by Proposition 2.2 (3), take $\Omega=\Lambda$,
 $A_{\lz}={\cal F}_{\lz}(E)$, $A={\cal F}(E)$. By Lemma 3.5 and Lemma 3.6,
 ${\cal F}_{\lz}(E)$ is contained in the $AF$ subalgebra ${\cal F}_{\lz}$
 of $A={\cal F}(E)$. Therefore ${\cal F}(E)$ is an $AF$ algebra. If moreover
 $G$ is compact, and so $\Gamma$ is discrete, we let
 $\mathbf{I}=\Gamma$, and $U_i=\{i\}$ for any
     $i\in\mathbf{I}=\Gamma$, then
     $\;C^*(E) \times_{\az^{\oz}}G=
     {\cal F}(E)$ as in Proposition 2.2. This completes the proof.\\

From the following theorem one can see how much difference between
the $AF$-embedding sufficient condition in Theorem 3.8 and the
$AF$-embedding necessary condition in Proposition 2.4 in our special case.\\

 {\sl{\bf Theorem 3.9.} Let $T$ be a finite
directed graph whose each weak connected component is a rooted
directed tree, $E$ be a finite directed graph constructed by
attaching $n_v$ ($0\le n_v <+\infty$) 1-loops to each vertex $v$
in $T$,  $G$ be a compact abelian group with the dual $\Gamma$,
$\oz: E^*\to \Gamma$ be a labelling map. Then the following are equivalent: \\
(1) $\oz: E^*\to \Gamma$ is almost proper.\\
(2) $\sum_{j=1}^k \oz_{\gamma_j}\not= 1_{\Gamma}$, for any loops
$\gamma_1$, $\gamma_2$, $\cdots$, $\gamma_k$ which are attached to
one path $\mu_0$ in $T$. }\\

 {\bf Proof:}(1)$\Rightarrow$(2):  obvious.\\
  (2)$\Rightarrow$(1): If $\oz$ is not almost proper, i.e. there is
 a compact subset $A$ of $\Gamma$ such that $\oz^{-1}(A)\backslash
 E^0$ isn't a finite set. Since $\Gamma$ is compact, $A$ is a
 finite set, and so there is an element $\sigma\in \Gamma$ such
 that $\oz^{-1}(\{\sigma\})\backslash
 E^0$ is a countable infinite set, denoted by $\{\nu_1, \nu_2,
 \cdots, \nu_n,\cdots\}$. For $i\ge 1$, let $\mu_i=\nu_i\cap T$,
 which consists of all the edges both in $\nu_i$ and in $T$, then
 $\mu_i$ ($i\ge 1$) are the paths in $T$ by the construction of $E$.
 Since $T$ is finite, there are only finite pathes in $T$. Without
 loss of generality, we may assume there is a path $\mu_0$
  in $T$ such that $\mu_i=\mu_0$ for all $i\ge 1$. Since $E^1$ is
  finite, by the similar discussion as in the proof of Proposition
  2.6, we may moreover assume that
  $\{l_i(e)\}_{i\ge 1}$ ($\forall e\in E^1$) is increasing,
    $\{|\nu_i|\}_{i\ge 1}$ is strictly increasing and
    $\lim_{i\to+\infty}|\nu_i|=+\infty$, where $l_i(e)$ is the
    appearing times of $e$ in $\nu_i$. Then by the construction of
    $E$ and that $\mu_j=\mu_0$ ($\forall j\ge 1$),
$\nu_{i+1}\backslash \nu_i$, which consists of all the edges in
$\nu_{i+1}$ but not in $\nu_i$, is a set of finite loops attached
to the path $\mu_0$. Therefore $\sum\limits_{\gamma\in
\nu_{i+1}\backslash \nu_i}\oz_{\gamma}\not=1_{\Gamma}.$ But, since
$\{l_i(e)\}_{i\ge 1}$ ($\forall e\in E^1$) is increasing,
$$\sum\limits_{\gamma\in \nu_{i+1}\backslash \nu_i}\oz_{\gamma}
=\oz_{\nu_{i+1}}-\oz_{\nu_{i}}=\sigma -\sigma=1_{\Gamma}.$$ This
is a contradiction, and completes the proof.

\vskip 1cm

{\bf Acknowledgements:}\\

This paper was written while the author visited the Mathematisches
Institut, Universitaet Muenster in the communication program of
DAAD and the Education Ministry of China. He is grateful to
Professor Joachim Cuntz for his hospitality, encouragement and
many important suggestions, and to Professor Siegfried Echterhoff
for his kind help. He is grateful to Professor George A. Elliott
for many important suggestions by email communication. The author
is  also pleasure to express his gratitude to DAAD and the
Education Ministry of China for their kind financial support. This
article is supported by National Natural Science Foundation of
China (10271090). XIAOCHUN FANG\\
Department of Applied Mathematics\\
Tongji University\\
Shanghai 200092, China\\
E-mail address:  xfang@mail.tongji.edu.cn

\end{document}